\documentclass[a4paper, 11pt]{amsart}

\usepackage{amsmath}
\usepackage{amsthm}
\usepackage{comment}

\usepackage[utf8]{inputenc}

\usepackage{xcolor}
\definecolor{green}{RGB}{0,127,0}
\definecolor{red}{RGB}{105,89,205}
\usepackage[colorlinks=true]{hyperref}

\theoremstyle{plain}
\newtheorem{thm}{Theorem}[section]

\theoremstyle{definition}

\newtheorem{ex}[thm]{Example}

\newtheorem{uw}[thm]{Remark}

\theoremstyle{remark}

\numberwithin{figure}{section}
\numberwithin{equation}{section}

\newcommand{\N}{\mathbb{N}}

\newcommand{\F}{\mathbb{F}}

\newcommand{\U}{\mathcal{U}}

\counterwithin{equation}{section}

\DeclareMathOperator{\GCD}{GCD}
\DeclareMathOperator{\ord}{ord}

\begin{document}

\title[Binomials and algebraic closure of  pseudofinite fields]{On binomials and algebraic closure of some pseudofinite fields}

\author{Jakub Gismatullin}

\address{Instytut Matematyczny Uniwersytetu Wroc{\l}awskiego, pl. Grunwaldzki 2/4, 50-384 Wroc{\l}aw, Poland \& Instytut Matematyczny Polskiej Akademii Nauk, ul. {\'S}niadeckich 8, 00-656 Warszawa, Poland}

\email{jakub.gismatullin@uwr.edu.pl}

\thanks{\noindent The first author is supported by the National Science Centre, Poland NCN grant no.  2017/27/B/ST1/01467.}

\author{Katarzyna Tarasek}

\address{Instytut Matematyczny Uniwersytetu Wroc{\l}awskiego, pl. Grunwaldzki 2/4, 50-384 Wroc{\l}aw, Poland}

\email{ja-takasia@wp.pl}

\date{\today}

\keywords{binomial, pseudofinite field, algebraic closure}
\subjclass[2010]{12L10, 12F05, 03C20}

\begin{abstract}
We give a criterion when a polynomial $x^n-g$ is irreducible over a pseudofinite field. As an application we give an explicit description of algebraic closure of some pseudofinite fields of zero characteristic.
\end{abstract}

\maketitle

	%\tableofcontents
	
	\section{Introduction}
	
	In this note we study binomials over pseudofinite fields. A \emph{binomial} over a field $K$ is a polynomial of the form $x^n-g$, where $n\in\N$ and $g\in K$. Theorem \ref{thm:3} below gives a criterion for irreducibility of a binomial over an ultraproduct of finite fields, that is \emph{pseudofinite field}. We construct in Example \ref{ex:f} a pseudofinite field $\F_{\bar{q}}$ and find $g\in \F_{\bar{q}}$, such that $x^n-g$ is irreducible over $\F_{\bar{q}}$ for all $n\in\N$. We then prove in Theorem \ref{thm:a}, that the algebraic closure $\widehat{\F_{\bar{q}}}$ is of the form $\F_{\bar{q}}\left( \sqrt[n]{g} : n\in\N\right)$.
	
	\section{Binomials and algebraic closure}\label{bn}
	
	Suppose $K$ is a field. A \emph{binomial} is a polynomial $w(x)\in K[x]$ of the form
	\[w(x)=x^n-g,\]
	for some $n\in\N$ and $g\in K$. 
	
	We refer the reader to \cite[Section 7, 7.9]{fa} and \cite{zoe} as a reference for basic notions around pseudofinite fields, that is ultraproducts of finite fields.
    
    Let us fix a nonprincipal ultrafilter $\U$ on $\N$. Let ${\mathcal F}=\{\F_{q_k}\}_{k\in\N}$ be a family consisting of some finite fields, where each $q_k=p_k^{t_k}$ is a power of a prime $p_k$, $t_k\in\N_{>0}$ and $\lim_{k\to\infty}q_k = \infty$.    Consider an ultraproduct 
    \[\F_{\bar{q}} = \prod_{k\in\N} \F_{q_k}/\U, \text{ where }\bar{q}=(q_k)_{k\in\N}\]
    of $\mathcal F$. Then $\F_{\bar{q}}$ is an infinite pseudofinite field. We study the following question:
	\begin{quote}
	    fix a pseudofinite field $\F_{\bar{q}} = \prod_{k\in\N} \F_{q_k}/\U$, is there $g\in \F_{\bar{q}}$ such that for every natural $n\in\N$, $x^n-g$ is irreducible over $\F_{\bar{q}}$?
	\end{quote}
	
   Let us recall a standard result on irreducible binomials from \cite[p. 425, Theorem 1.6]{karp}: 
   \begin{quote}
    Fix a field $K$, $\N\ni n \ge 1$ and $g \in K$. Then $x^n-g$ is irreducible over $K$ if and only if $(1)$ and $(2)$ hold:
    \begin{enumerate}
        \item $g\notin K^p$, for each prime $p$, which divides $n$,
        \item if $4|n$, then $g \notin -4K^4=\{-4x^4 : x\in K\}$. 
    \end{enumerate}
    \end{quote}
     This result for finite field has the following form:

\begin{thm}\cite[Theorem 3.75, p. 125]{lidl} \label{lidl}
    Let $q=p^m$ be a power of a prime $p$ and fix $w(x)= x^n - g$, for $g\in \F_q$, $n \ge 2$. Let $e$ be the order of $g$ in the multiplicative group $\F_q^\times=\F_q\setminus\{0\}$. Then $w$ is irreducible over $\F_q$ if and only if the following three conditions hold:
    \begin{enumerate}
        \item[1)] $\GCD((q-1)/e,n)=1$,
        \item[2)] every prime divisor of $n$ divides also $e$,
        \item[3)] if $4 | n$, then also $4 | q- 1$.
    \end{enumerate}
\end{thm}

\begin{uw}\label{uwaga}
Observe that $\F_q^\times$ is a cyclic group, that is $\F_q^\times$ is generated by a single element $g_q\in \F_q^\times$. When $x^n - g_q$ is irreducible over $\F_q$? Since $g_q$ has order $q-1$ in $\F_q^\times$, the condition 1) from theorem \ref{lidl} is always true. Therefore $x^n - g_q$ is irreducible over $\F_q$ if and only if  
\[ \text{every prime divisor of $n$ divides $q-1$ and if $4|n$, then $4|q-1$.} \tag{$\diamondsuit$} \]
For example, if $q=13$, then each prime divisor of $n$ must divide $q-1=12=2^2\cdot 3$. Hence $n$ must be of the form $2^m\cdot 3^k$, $m,k\in \N$. Every binomial $x^{2^m\cdot 3^k}-2$ is irreducible over $\F_{13}$ (here $g=2$ is a generator of $\F_{13}^\times$).
\end{uw}

\if
Let us give an illustration of \ref{lidl} for the following question
    \begin{quote}
        Fix a power $q=p^m$. For which $n\in\N$, there is $a\in\F_q$ such that $x^n-a$ is irreducible over $\F_q$.
    \end{quote}
    
    Theorem \ref{lidl} and Remark \ref{uwaga} imply the following partial answer.
    \[\begin{array}{| c|c|} \hline
 q & n   \\ \hline
 3 & 2  \\ \hline
 5 & 2^m   \\ \hline
 7 & 2^m 3^k, 4\not|n  \\ \hline
 9 & 2^m \\ \hline
 11 & 2^m 5^k, 4\not|n \\ \hline
 13 & 2^m 3^k \\ \hline
 17 & 2^m  \\ \hline
 19 & 2^m 3^k, 4\not|n \\ \hline
 23 & 2^m 11^k, 4\not|n \\ \hline
 25 & 2^m 3^k \\ \hline
 27 & 2^m 13^k, 4\not|n \\ \hline
 29 & 2^m 7^k \\ \hline
 31 & 2^m 3^k 5^w, 4\not|n \\ \hline
 37 & 2^m 3^k \\ \hline
 41 & 2^m 5^k \\ \hline
 43 & 2^m 3^k 7^w, 4\not|n \\ \hline
 47 & 2^m 23^k, 4\not|n \\ \hline
 49 & 2^m 3^k \\ \hline
\end{array}.\]
\fi

        We need the following elementary remark on definability of irreducibility of polynomials of a fixed degree. The proof is standard, and can be found e.g. in \cite[p. 136]{fa}.
        
        \begin{uw}\label{uww}
        For each natural $n$ there exists a first order formula $\varphi_n(x_1,\ldots,x_n)$ in the language of fields $L_F=\{+,\cdot,0,1\}$ such that, for every field $K$ and every collection of elements $a_0,\ldots,a_{n-1} \in K$
        \begin{quote}
        polynomial $x^n+a_{n-1}x^{n-1}+\ldots+a_1 x+a_0$ is irreducible over $K$ if and only if $\varphi_n(a_0,\ldots,a_{n-1})$ holds in $K$. 
        \end{quote}
        \end{uw}

        %Wykorzystam ponizszy fakt opisujacy  podgrupy grupy multyplikatywnej $\F^\times_q\cong \Z_{q-1}$ postaci         \[{\F^\times_q}^{p}=\{f^{p} : f\in \F^\times_q\},\]         gdzie $p$ jest liczba pierwsza, a $q={p'}^t$ jest potega jakiejs liczby pierwszej $p'$. Poniewaz $\F^\times_q$ jest grupa cykliczna izomorficzna z $\Z_{q-1}$, prawdziwy jest nastepujacy dobrze znany fakt .
        
        Theorem \ref{thm:3} below is our criterion for the irreducibility of $x^n-g$ over a pseudofinite field $\F_{\bar{q}}$. This criterion is a straightforward application of the Łoś ultraproduct theorem.
        
          \begin{thm} \label{thm:3}
            Let ${\mathcal F}=\{\F_{q_k}\}_{k\in\N}$ be a family of some finite fields and $\lim_{k\to\infty}q_k = \infty$. Take an ultraproduct $\F_{\bar{q}} = \prod_{k\in\N} \F_{q_k}/\U$. Fix $n\in\N$ and let $\{\text{prime divisors of }n\}=\{p_1,\ldots,p_r\}$.
            Then the following conditions are equivalent:
            \begin{enumerate}
                \item[(1)] $x^n-g$ is irreducible over $\F_{\bar{q}}$, where \[g=([g_k]_{k\in\N})_{\U}\in \F_{\bar{q}},\text{ where each }g_k\text{ is a generator of }\F_{q_k}^\times. \tag{$\clubsuit$}\]
                \item[(2)] there exists $g\in \F_{\bar{q}}$ such that $x^n-g$ is irreducible over $\F_{\bar{q}}$,
                \item[(3)] for $\U$-almost all $k\in\N$:
                \[p_1,p_2,\ldots,p_r\text{ divide }q_k-1\text{ and if }4|n\text{ then }4|q_k-1.\]
            \end{enumerate}
        \end{thm}
        
        By ``for $\U$-almost all $k\in\N$ the condition $\psi(k)$ holds'' we mean that $\{k\in\N : \text{  holds }\psi(k)\}\in\U$.
        
        \begin{proof}
        $(1)\Rightarrow(2)$ is clear. $(2)\Rightarrow(3)$  Suppose that $g=([g_k]_{k\in\N})_{\U}\in \F_{\bar{q}}$, where $g_k\in \F_{q_k}$ and $x^n-g$ is irreducible over $\F_{\bar{q}}$. By Remark \ref{uww} and Łoś's theorem, $x^n-g_k$ is irreducible over $\F_{q_k}$, for $\U$-almost all $k\in\N$. Theorem \ref{lidl} implies that: $p_1,p_2,\ldots,p_r$ divide $e=\ord(g_k)$, which is a divisor of $q_k-1$, and if $4|n$, then $4|q_k-1$ (for $\U$-almost all $k\in\N$). This is exactly $(2)$.
        
        $(3)\Rightarrow(1)$ Suppose $(3)$ is true. Let  $g=([g_k]_{k\in\N})_{\U}\in \F$ be as in $(\clubsuit)$. In order to obtain (1) it is enough (by Remark \ref{uww}) to prove that 
        \[\text{for $\U$-almost all $k\in\N$, the binomial $x^n-g_k$ is irreducible over  $\F_{q_k}$}.\] This is exactly (3), by $(\diamondsuit)$ in \ref{uwaga}.
%        Clearly $(\clubsuit)$ is equivalent to
%        \[\text{for $\U$-almost all $k\in\N$, }g_k\notin \F_{q_k}^{p_1}\cup\ldots\cup \F_{q_k}^{p_r}. \tag{$\spadesuit$}\]
%        
%        \begin{enumerate}
%        \item If $n$ is coprime to $q-1$, then $\F^\times_q={\F^\times_q}^{n}=\{f^n : f\in \F^\times_q\}$.
%        \item If $n$ divides $q-1$, then $|{\F^\times_q}^{p}| = \frac{q-1}{p}$ and for each $f\in \F^\times_q$
%            \[f\in {\F^\times_q}^{p} \Longleftrightarrow f^{\frac{q-1}{p}}=1.\]
%        \end{enumerate}    
        \end{proof}
        
        As an application we construct a pseudofinite field $\F_{\bar{q}}$ of characteristic zero and $g\in \F_{\bar{q}}$, such that $x^n-g$ is irreducible for every $n\in\N$.
        
        \begin{ex}\label{ex:f}
        Let $\{2=p_1,p_2,\ldots\}$ be the set of all prime numbers. Take a sequence $\{q_k : k\in\N\}$ of prime powers, such that
        \[4,p_1,p_2,\ldots,p_k \text{ divide } q_k-1\text{ for all natural }k\in\N. \tag{$\spadesuit$}\]
        For example one can take
        \[q_k = p_{k+1}^{(p_1-1)(p_2-1)\cdot\ldots\cdot(p_k-1)},\] as then $p_i|q_k-1={r_k}^{p_i-1}-1$ by the Fermat's little theorem.
        Then Theorem \ref{thm:3} can be applied to $\F_{\bar{q}} = \prod_{k\in\N} \F_{q_k}/\U$, to get that $x^n-g$ is irreducible over $\F_{\bar{q}}$, for all natural $n\in\N$.
        \end{ex}
        
        Another application of Theorem \ref{thm:3} and Example \ref{ex:f} is an explicit description of algebraic closure of $\F_{\bar{q}}$ from Example \ref{ex:f}. 
        
        \begin{thm}\label{thm:a}
        Suppose $\F_{\bar{q}}$ is a pseudofinite field such that for some $g\in \F_{\bar{q}}$, for every natural $n\in\N$, $x^n-g$ is irreducible over $\F_{\bar{q}}$. Then its algebraic closure $\widehat{\F_{\bar{q}}}$ is generated by $\left\{ \sqrt[n]{g} : n\in\N\right\}$ over $\F_{\bar{q}}$:
        \[\widehat{\F_{\bar{q}}} = \F_{\bar{q}}\left( \sqrt[n]{g} : n\in\N\right),\]
        where $\sqrt[n]{g}$ is a root of $x^n-g$.
        \end{thm}
        \begin{proof}
        We use a well know fact, that every pseudofinite field has a unique algebraic extension of each natural degree \cite{zoe}.
        Suppose that $a\in \widehat{\F_{\bar{q}}}$ has degree $n$ over $\F_{\bar{q}}$, that is $[\F_{\bar{q}}(a):\F_{\bar{q}}]=n$. Then since also $\left[\F_{\bar{q}}\left(\sqrt[n]{g}\right):\F_{\bar{q}}\right]=n$, we have that $\F_{\bar{q}}(a) = \F_{\bar{q}}\left(\sqrt[n]{g}\right)\subseteq \F_{\bar{q}}\left( \sqrt[n]{g} : n\in\N\right)$.
        \end{proof}
        
        A pefect field $K$ is called \emph{quasi-finite}, if $K$ has a unique (necessarily cyclic) extension of degree $n$ for each $n\in\N$. Theorem \ref{thm:a} is true for quasi-finite fields, with the same proof.

\bibliographystyle{alpha}
\bibliography{literatura}	

\end{document}